\documentclass[11pt]{article}
\newtheorem{theorem}{Theorem}
\newtheorem{lemma}{Lemma}
\usepackage{latexsym}
\begin{document}
\date{}
\title{A combinatorial determinant}
\author{Herbert S. Wilf\\\normalsize Department of Mathematics, University
of Pennsylvania\\\normalsize Philadelphia, PA 19104-6395}
\maketitle
\begin{abstract}
A theorem of Mina evaluates the determinant of a matrix with entries
$D^j(f(x)^i)$. We note the important special case where the matrix entries
are evaluated at $x=0$ and give a simple proof of it, and some
applications. We then give a short proof of the general case.
\end{abstract}\vspace{.1in}
\begin{center}
October 25, 1998
\end{center}
\vspace{.1 in}
An old theorem of Mina \cite{Mi}, which ``deserves to be better known''
\cite{vdp}, states that
\begin{equation}
\label{eq:mena}
\det{\left\{\left(\frac{d^j}{dx^j}f(x)^i\right)_{i,j=0}^{n-1}\right\}=
1!2!\dots (n-1)!f'(x)^{{n(n-1)/2}}.         }
\end{equation}
A proof of Mina's theorem can be found, for instance, in \cite{WS}.
We will first give a short proof of the special case in which both sides
are evaluated at $x=0$, with some applications, and then give a short proof
of the general case. The special case shows an interesting structure owing
to the fact that all matrices of that form can be simultaneously
triangularized by multiplying them by a certain universal triangular matrix.
\section{Coefficients of powers of power series}
\begin{theorem}
\label{th:main}
Let $f=1+a_1x+a_2x^2+\dots$ be a formal power series, and define a matrix
$c$ by\footnote{ \ ``$[x^k]g$'' means the coefficient of $x^k$ in the
series $g$.}
\begin{equation}
\label{eq:defmat}
c_{i,j}=[x^j]f^i\qquad (i,j\ge 0).
\end{equation}
Then
\begin{equation}
\label{eq:main}
\det{\left((c_{i,j})_{i,j=0}^n\right)}=a_1^{n(n+1)/2}\qquad
(n=0,1,2,\dots ).
\end{equation}
\end{theorem}

To prove this, define a matrix $b$ by
\begin{equation}
\label{eq:bdef}
b_{i,j}=(-1)^{i+j}{i\choose j},\qquad (i,j\ge 0).
\end{equation}
 Then we claim that $bc$ is upper triangular
with powers of $a_1$ on its diagonal. Indeed we have
\begin{eqnarray*}
\sum_jb_{i,j}c_{j,k}&=&\sum_j(-1)^{i+j}{i\choose
j}[x^k]f^j=(-1)^i[x^k]\sum_j(-1)^j{i\choose j}f^j\\
&=&(-1)^i[x^k](1-f)^i=[x^k](a_1x+a_2x^2+\dots )^i=\cases{0,&if $k<i$;\cr
a_1^i,&if $k=i$,\cr}
\end{eqnarray*}
as claimed. Since $bc$ is this upper triangular matrix, and $\det{b}=1$,
the determinant of $c$ is $\prod_{i=0}^na_1^i=a_1^{n(n+1)/2}$. $\Box$

Ed Bender has noted that the hypothesis $a_0=1$ can be removed. Indeed if
$a_0\neq 0$, apply the result to $f/a_0$ and discover that the theorem is
unchanged. If $a_0=0$ the result follows by continuity.

In order to gain an extra free parameter in the identities that are to
follow, as well as to introduce the idea of the proof of Mina's theorem in
general form, we'll restate Theorem \ref{th:main} in terms of the $z$th
power of $f$.
\begin{theorem}
\label{th:main2}
Let $f=1+a_1x+a_2x^2+\dots$ be a formal power series, let $z$ be a complex
number, and define a matrix $c$ by
\begin{equation}
\label{eq:defmat2}
c_{i,j}=[x^j]f^{zi}\qquad (i,j\ge 0).
\end{equation}
Then
\begin{equation}
\label{eq:main2}
\det{\left((c_{i,j})_{i,j=0}^n\right)}=(za_1)^{n(n+1)/2}\qquad
(n=0,1,2,\dots ).
\end{equation}
\end{theorem}

\section{Some examples}
\begin{enumerate}
\item This investigation began when I was looking at the infinite matrix
whose $(i,j)$ entry is the number of representations of the integer $j$ as
a sum of $i$ squares of nonnegative integers $(i,j=0,1,2,\dots)$, and
noticed that its determinant is 1. This matrix begins as
\[\left(\begin{array}{rrrrrrrrr}
1& 0& 0& 0& 0& 0& 0& 0&\dots\\ 1& 1& 0& 0& 1& 0& 0& 0&\dots\\ 1& 2& 1& 0&
2& 2& 0& 0&\dots\\
   1& 3& 3& 1& 3& 6& 3& 0&\dots\\ 1& 4& 6& 4& 5& 12& 12& 4&\dots\\ 1& 5&
10& 10& 10& 21& 30& 20&\dots\\
 1& 6& 15& 20& 21& 36& 61& 60&\dots\\ 1& 7& 21& 35& 42& 63& 112& 141&\dots\\
\vdots&\vdots&\vdots&\vdots&\vdots&\vdots&\vdots&\vdots&\ddots
\end{array}\right)\]
This is an array of the type considered in (\ref{eq:defmat}) above, where
$f=1+x+x^4+x^9+x^{16}+\dots$. Hence by Theorem \ref{th:main}, the
determinant of every upper-left $n\times n$ section of this infinite matrix
is 1. The same will be true if ``squares'' is replaced by ``cubes'' or any
higher power, or for that matter by any increasing sequence $\{0,1,\dots
\}$ at all!
\item Take $f=1+x$ in Theorem \ref{th:main2} to discover that
\[\det{\left({zi\choose j}\right)_{i,j=0}^n}=z^{{n+1\choose 2}}\qquad
(n=0,1,2,\dots ).\]
\item With $f=(e^x-1)/x$ in Theorem \ref{th:main2} we find a determinant
that involves the Stirling numbers of the second kind,
\[\det{\left(\frac{(zi)!}{(zi+j)!}{zi+j\brace
zi}\right)_{i,j=0}^n}=\left(\frac{z}{2}\right)^{n(n+1)/2}\qquad
(n=0,1,2,\dots ).\]
\item With $f=\log{(1+x)}/x$ in Theorem \ref{th:main2} we evaluate one that
contains the Stirling numbers of the first kind,
\[\det{\left(\frac{(zi)!}{(zi+j)!}{zi+j\brack
zi}\right)_{i,j=0}^n}=\left(\frac{z}{2}\right)^{n(n+1)/2}\qquad
(n=0,1,2,\dots ).\]
\item Now let $f=(1-\sqrt{1-4x})/(2x)=1+x+\dots$ in Theorem \ref{th:main2}.
If we use the fact that
\[ \left(\frac{1-\sqrt{1-4x}}{2x}\right)^k=\sum_{n\ge
0}\frac{k(2j+k-1)!}{j!(k+j)!}x^n,\]
then the resulting determinantal identity can be put in the form
\[\det{\left(\frac{(2j+zk-1)!}{(zk+j)!}\right)_{j,k=1}^n}=z^{{n\choose
2}}\,1!\,2!\,3!\,\dots (n-1)!\qquad (n=1,2,\dots ).\]
\end{enumerate}

\section{The additive structure of this case}

This case of Mina's theorem, in which we evaluate at $x=0$, has an
interesting structure because, if the series involved are normalized to
constant term 1, then the matrices are all trangularized by the same matrix
$b$ of (\ref{eq:bdef}).

Usually, determinants don't relate well to addition of matrices. That is,
if $\det{U}$, $\det{V}$ can be evaluated in simple closed form, there is no
reason to suppose that the same is true of $\det{(U+V)}$. But the
determinants that we are now studying act nicely under matrix addition.
Thus, if $u_{i,j}=[x^i]f^j$, and $v_{i,j}=[x^i]g^j$, then $b(u+v)$ is
triangular, with diagonal entries $f'(0)^i+g'(0)^i$ $(i=0,1,2,\dots)$, and so
\[\det{(u_{i,j}+v_{i,j})_{i,j=0}^n}=\prod_{i=0}^n(f'(0)^i+g'(0)^i).\]

For example, \[\det{\left({ri\choose j}+{si\choose
j}\right)_{i,j=0}^n}=\prod_{i=0}^n(r^i+s^i),\]
and more generally if $S$ is any set of numbers, and $c$ is a function on
$S$, then
\[\det{\left(\sum_{r\in S}c(r){ri\choose
j}\right)_{i,j=0}^n}=\prod_{j=0}^n\left(\sum_{r\in S}c(r)r^j\right).\]
An interesting special case is obtained by taking $S$ to be a set of $m$
equally spaced points in an interval, say $(0,1)$, taking the $c(r)$'s all
equal to $1/m$, and taking the limit as $m\to\infty$. The result is that
\[\det{\left(\int_0^1f(x){xi\choose
j}dx\right)_{i,j=0}^n}=\prod_{j=0}^n\mu_j(f),\]
where the $\mu_i(f)=\int_0^1x^if(x)dx$ are the moments of $f$. For example,
\[\det{\left(\int_0^1{xi\choose j}dx\right)_{i,j=0}^n}=\frac{1}{(n+1)!}.\]

\section{Proof of Mina's theorem}
Now we prove Mina's theorem by proving the following small generalization
of it.
\begin{equation}
\label{eq:menag}
\det{\left\{\left(\frac{d^j}{dx^j}f(x)^{zx_i}\right)_{i,j=0}^{n}\right\}=f(x
)^{z\sum_0^{n}(x_i-i)}\left(zf'(x)\right)^{n(n+1)/2}\prod_{0\le i<j\le
n}(x_j-x_i).         }
\end{equation}
This will follow easily from the following observation.
\begin{lemma}
\label{lm:polys1}
 Let $\{p_j(x)\}$ $(j=0,1,2,...)$ be any sequence of polynomials such that
the $j$th one is of degree $j$, for each $j\ge 0$. Then the determinant of
the matrix $Q=\{p_j(x_i)\}_{i,j=0}^{n}$ is equal to the product of the
highest coefficients of $p_0,p_1,\dots ,p_{n}$ times the discriminant of
the $x_i$'s.
\end{lemma}
Indeed, if $A$ is the lower triangular matrix of coefficients of the
$p_j$'s and $V$ is the Vandermonde matrix of the $x_i$'s then $Q=AV$. $\Box$

Now let $\alpha_{i,j}=D^j(f(x)^{zx_i})$, where $D=d/\kern-2pt dx$. We claim
that, \textit{for each $j=0,1,\dots,$ $\alpha_{i,j}$ is of the form
$f^{zx_i-j}$ times a polynomial in $z$ of exact degree $j$, whose
coefficients are functions of $x$ and whose highest coefficient is
$(x_if')^j$}. In view of Lemma \ref{lm:polys1}, this claim will prove the
theorem. However this claim is trivial to prove by induction on $j$. $\Box$

Theorem \ref{th:main2} is the special case of (\ref{eq:menag}) in which
$x=0$, $z=1$, and the polynomials $p_j(x)$ are given by
\[p_j(x)=\sum_{\ell=0}^j\frac{x^{\ell}}{\ell!}[t^j]
\left(\log{y(t)}^{\ell}\right),\]
with highest coefficients $a_1^j/j!$.

\section{Remarks}
\begin{enumerate}
\item The family of matrices of the form (\ref{eq:defmat}) was previously
investigated by
Doron Zeilberger \cite{DZ}, from the point of view of constant term
identities.
\item My thanks go to Brendan McKay for some e-mail exchanges that were
helpful in clarifying the ideas here.
\item By comparing (\ref{eq:mena}) and (\ref{eq:menag}) we see that the
absence of $f(x)$ on the right side of the former is because of the ``lucky
coincidence'' that in the former, $x_0+x_1+\dots +x_{n}=n(n+1)/2$.
\end{enumerate}


\begin{thebibliography}{aaaa}
\bibitem{Mi} L. Mina, \textit{Formole generali delle derivate successive d'una
funzione espresse mediante quelle della sua inverse}, Giornale
di Mat. xliii (1904), 196-212.
\bibitem{WS} Volker Strehl and Herbert S. Wilf, Five surprisingly simple
complexities, Symbolic computation in combinatorics $\Delta_1$ (Ithaca, NY,
1993), \textit{J. Symbolic Comput.} \textbf{20} (1995), no. 5-6, 725--729.
\bibitem{vdp} A. J. van der Poorten, Some determinants that should be
better known, {\it J. Austral. Math. Soc.} {\bf 21} (Series A) (1976),
278-288.
\bibitem{DZ} Doron Zeilberger, A constant term identity featuring the
ubiquitous (and mysterious)
Andrews-Mills-Robbins-Rumsey numbers $\{1,2,7,42,429, ...\}$,
\textit{J. Combinatorial Theory (Ser. A)} \textbf{66}  (1994), 17--27.
\end{thebibliography}
\end{document}